\documentclass{article}

\usepackage{PRIMEarxiv}

\usepackage[utf8]{inputenc} 
\usepackage[T1]{fontenc}    
\usepackage{hyperref}       
\usepackage{url}            
\usepackage{booktabs}       
\usepackage{amsfonts}       
\usepackage{nicefrac}       
\usepackage{microtype}      
\usepackage{lipsum}
\usepackage{fancyhdr}       
\usepackage{graphicx}       
\usepackage{subcaption}
\graphicspath{{media/}}     
\title{Restoring Accessibility During Urban Rail Disruptions via Bus Network Redesign
}
\usepackage{authblk}

\author[1]{Zihao Guo}
\author[1]{Andrea Araldo}
\author[1]{Mounîm A. El Yacoubi}

\affil[1]{Samovar, Télécom SudParis, Institut Polytechinique de Paris, Évry-Courcouronnes, France\\
\texttt{zihao-eric.guo@ip-paris.fr}, 
\texttt{\{andrea.araldo, mounim.el\_yacoubi\}@telecom-sudparis.eu}}

\newcommand{\keepcomment}{1} 
\ifnum\keepcomment=0
	\usepackage[colorinlistoftodos,textsize=scriptsize]{todonotes} 
    \newcommand{\stkout}[1]{\ifmmode\text{\sout{\ensuremath{#1}}}\else\sout{#1}\fi}
    
\else
	
	\usepackage[disable]{todonotes} 
\fi

\usepackage{amsmath}
\usepackage{algorithm}
\usepackage{algorithmic}
\usepackage{calrsfs}
\DeclareMathAlphabet{\pazocal}{OMS}{zplm}{m}{n}
\usepackage{tabularx}
\usepackage{subcaption}
\usepackage{moreverb,xurl}

\begin{document}
\maketitle

\begin{abstract}
In broad terms, accessibility measures opportunities reachable (such as shops, residents, etc.) within a given time frame. Urban Rail Transit (URT) plays a crucial role in providing accessibility, but it is susceptible to disruptions. In city centers with dense public transport (PT) networks, travelers can often find alternative lines. However, in suburbs where PT is sparse, disruptions have a more significant impact on accessibility.
The traditional approach consists in deploying bridge and replacement buses to mitigate URT disruptions without specific care to accessibility. Yet, the question arises: is this approach the most effective way to restore accessibility? To the best of our knowledge, our paper is the first to propose a bus re-routing method with the objective of restoring accessibility during URT disruptions.
We formulate an integer program and develop a two-stage heuristic algorithm to maximize restored accessibility. The efficacy of our method is always the present assessed in Évry-Courcouronnes and Choisy-le-Roi, France. The results show that, compared to conventional replacement methods, our strategy improves accessibility in particular in the areas most affected by the disruption. Such results are observed even when no additional vehicles are deployed, and at the same time, achieving a reduction in the kilometers traveled. 
Despite it is well understood that accessibility is the most relevant benefit a transportation system can produce, this aspect is reflected by the traditional approaches in remediation to disruption. With this work, we show instead how to make accessibility the main guiding principle in remediation.
\end{abstract}


\section{Introduction}
Urban Rail Transit (URT) plays a pivotal role in metropolitan areas by alleviating suburban-to-downtown congestion, offering high energy efficiency, and convenient service across various weather conditions~\cite{ieda1995commuter,tan2020evacuating}. However, URT systems face operational challenges to maintain high reliability~\cite{pender2012planning}. Notable incidents, such as the July 2017 fire at New York City's 145th Street Station caused by track debris, have led to significant service disruptions for more than 2 hours during morning peak~\cite{barone2017abcdfire}. Even the highly efficient Hong Kong Mass Transportation Railway (MTR), boasting a 99.9\% on-time performance, encounters approximately 250 disruptions annually~\cite{zhang2020metro}. In Paris, the RER B line typically serves nearly 200,000 passengers on a weekday. However, for maintenance work, long segments of the line have been closed in recent years, over multiple hours. Replacement buses operating to remediate this closure have proven to be insufficient, forcing many commuters to seek alternative routes or reduce their travel~\cite{actu2024travaux}.

In city centers with dense public transport (PT) networks, passengers are typically able to find alternative lines. With sparse public transport, nonetheless, the effect of URT disruption will be much worse. Due to the lack of alternatives, passengers in these areas might be left with no mobility options and have severe difficulties in performing their activities. To take these aspects into account, the usual level of service metrics (average waiting times, travel times, etc.) are not sufficient. What is more relevant to passengers is instead \emph{accessibility}. Stated in simple terms, accessibility indicators measure how many opportunities per hour can be reached from a certain location~\cite{miller2020measuring}. Opportunities can be workplaces, schools, healthcare facilities or any other location that must be reached to perform activities relevant for the users. 
Nonetheless, remediation strategies applied during disruptions do not generally take accessibility into account. 
Indeed, conventional remediation strategies consist of providing substitute buses and primarily focus on minimizing operational metrics while overlooking the crucial factor of area accessibility. This represents a significant gap in current practice, as maintaining accessibility is fundamental to preserving residents' quality of life and social participation during disruptions~\cite{liu2024measuring}.

To fill this gap, we propose a bus network redesign method to restore accessibility during URT disruptions. Our key innovation lies in making accessibility the primary optimization objective when redesigning bus networks during rail disruptions, rather than focusing solely on conventional operational metrics. In order to serve the impacted passengers, our approach decides which lines to extend, calculates the routes of these extensions, and re-optimizes the allocation of buses across lines. 
We show (and explain) the superiority of our approach with respect to usual practice, which consists in just replacing the interrupted URT line with a ``replacement bus'', which just stops at the original disrupted stops. Due to the problem's complexity, we further develop a two-stage heuristic optimization method that generates solutions within seconds. We validate our approach using two case studies of rail disruptions in two French towns in the suburbs of the Paris region. Results demonstrate that our method recovers the accessibility loss during disruptions much more efficiently than conventional replacement buses, restoring more accessibility with less driving distance, even without adding buses, just by appropriately extending existing lines and reallocating the existing fleet of buses.

We pinpoint that we are not tackling in this paper emergency situations, as in those cases, the priority would be to save lives rather than ensuring accessibility. We consider instead disruptions occurring on a longer timescale (hours or days, for instance - see real examples in~\cite{christoforou2016managing,20minutes2024}), which could be planned or unplanned and that have a big impact on the everyday life of commuters. With such timescales, it is reasonable to assume that users will have the time to be notified of the disruption and to recompute their journeys. We assume users will use a journey planner (in the form of a smartphone application), which can tell them the optimal journey within our re-arranged PT structure.

\section{Related Work}

Due to its inherent complexity, Urban Rail Transit (URT) disruption management has been largely studied in the literature. Shalaby et al.~\cite{shalaby2021rail} reviewed recovery models and algorithms for real-time railway disruption management, which involves timetable adjustment, substitute bus scheduling, and crew rescheduling.

\subsection{Work Related to Bus Bridging}
In the case where a disruption interrupts a rail line "in the middle" by keeping the two extreme segments isolated from each other, strategies often consist in deploying \emph{bridge buses} in order to reconnect those two segments, as in~\cite{kepaptsoglou2009transit}. The design of shuttle bus is a critical issue.
Wang et al.~\cite{wang2016feeder} introduced a model for bus dispatching and route design during emergencies. Chen and An~\cite{chen2021integrated} developed an integrated optimization model to address the issues of bus timetabling and bridging route design, six distinct bus bridging routes were created, each passing via a different disrupted station, to give customers at those stations a choice of routes to choose from.
and Xu et al.~\cite{xu2021optimizing} created a robust optimization model considering uncertain disruption duration, while Gu et al.~\cite{gu2018plan} designed express routes with stop-skipping capabilities for high-flow stations in order to ensure the flexibility of bus bridging services during disruptions.
Deng et al.~\cite{deng2018design} further expanded on this basis, considering various bridging types, namely different combinations of bus bridging and urban rail transit (including bus bridging followed by rail transit, rail transit followed by bus bridging, and rail transit-bus bridging-rail transit), and established a model with station capacity constraints. However, these conventional studies have not addressed the challenges of distance between spare bus sources and demand points. This can lead to bridging buses traveling long distances to perform tasks due to a lack of nearby vacant public transportation. Our method creates nearby bus sources by extending regular bus lines, thus reducing this deadhead distance.
Recently, some studies~\cite{wang2022real, feng2024bus} explored the use of in-service buses (buses that are already in operation on regular lines) for disruption management. The limit of bus bridging is that it cannot be applied when an entire route is disrupted.

\subsection{Work Related to Replacement Buses}
When an entire rail line is disrupted, an approach commonly employed in practice~\cite{actu2024travaux} is to deploy \textbf{replacement buses (also called dedicated shuttle fleets)} to travel across all stops of the disrupted rail line~\cite{pender2012proactive, van2016shuttle}.
Previous literature mainly focused on minimizing transit operating time, passenger discomfort, and system cost. In the same sense, Wang et al.~\cite{wang2016feeder} created a model for bus routes optimization for overall time driving minimization, whereas Cadarso et al.~\cite{cadarso2013recovery} designed an optimization model of bus bridging
schedules and vehicle rescheduling from a comprehensive point of view, considering minimization of time of recovery, passenger discomfort, and system costs. Luo and Xu~\cite{luo2021railway} incorporated the stochastic nature of passenger demand and existing rail and bus routes' backup capacities into a stochastic programming model to minimize expected unmet passenger requirements. To cope with uncertain and heterogeneous bus traveling times, Liang et al.~\cite{liang2019robust} developed a passenger flow and operational cost optimization model. 

\subsection{Accessibility Gap in Existing URT Disruption Literature}
Despite the amount and high quality of previous work on Urban Rail Transit disruption remediation, all previous studies have overlooked the most crucial factor of a transportation system, i.e., \emph{accessibility}~\cite{miller2020measuring}. This aspect is particularly important in suburban regions, where accessibility often tightly depends on one or few rail lines. Due to the scarcity of public transport alternatives, a disruption of a rail line can severely impact the accessibility of suburbs, preventing residents from performing their daily activities.
\\
While accessibility is a well-established concept in transportation planning, our investigation confirms that URT disruption management research has not explicitly incorporated it as an optimization objective. Liu et al.~\cite{liu2024measuring} define accessibility as ``the potential of opportunities for interaction'' and recognize that transit systems face high risk of degradation from external disruptions, affecting their ability to deliver reliable accessibility. Their work introduces measures to detect disruptions' impacts, but focuses on measurement rather than intervention. Similarly, Nalin et al.~\cite{nalin2025evaluation} emphasize the necessity of assessing the actual performance of public transportation rather than relying on schedules when evaluating accessibility, yet do not propose recovery strategies. 
This critical gap in the literature, the absence of optimization methods to restore accessibility through network redesign, constitutes our paper's core contribution.

\subsection{Our Contribution}

For this reason, to the best of our knowledge, we are the first to propose a bus network redesign method with the explicit objective of optimizing accessibility restoration during Urban Rail Transit (URT) disruptions, particularly suited for suburbs with scarce public transport resources.

\section{Methodology}
In this section, we first formulate the problem of network redesign to restore accessibility in case of disruption of URT. Secondly, we develop an integer programming model to obtain the optimal bus fleet allocation.
\begin{figure*}[!t]
    \centering
    \includegraphics[width=1\textwidth]{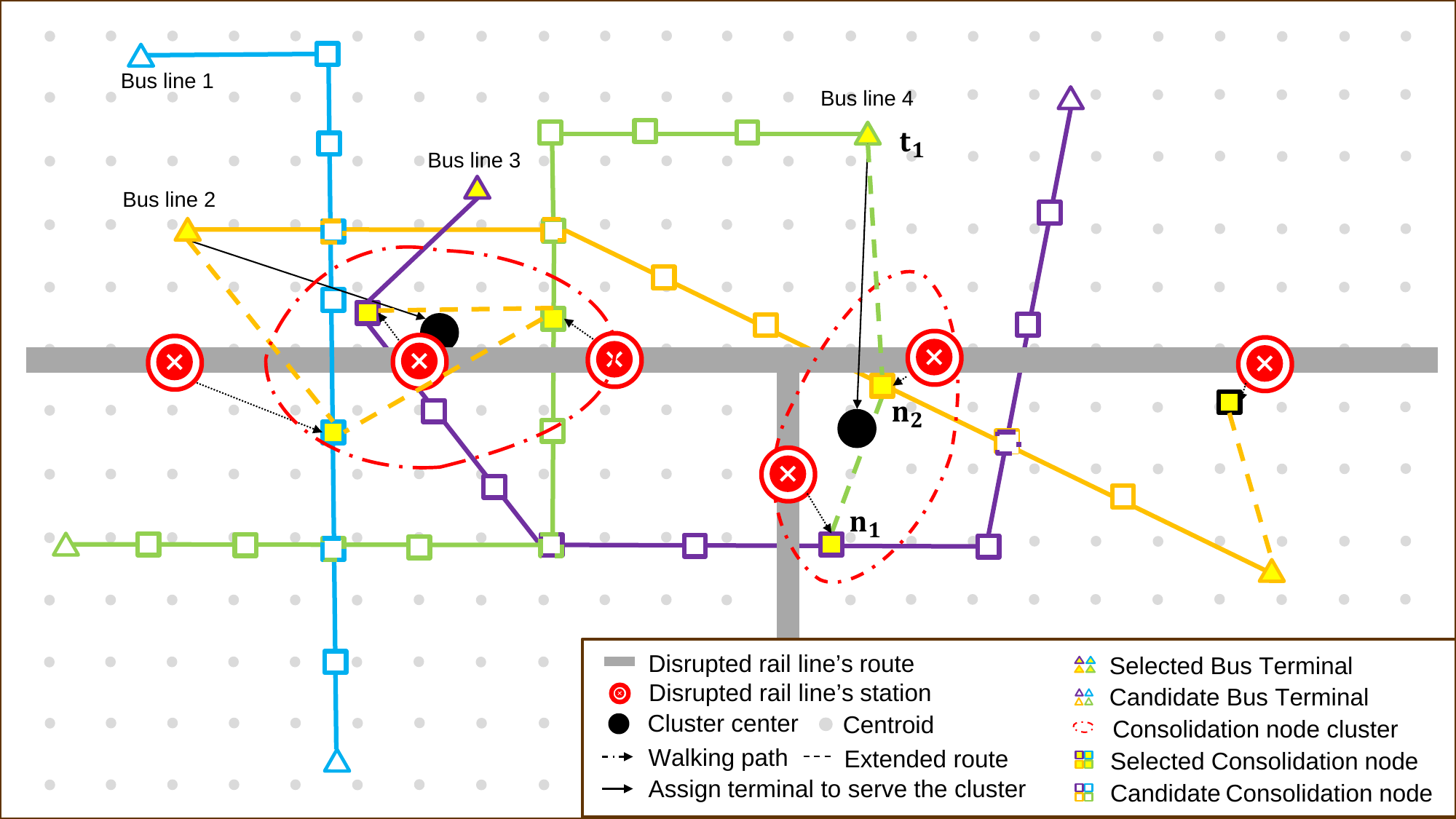}
    \caption{Network representation}
    \label{fig:PT_repre}
\end{figure*}

\subsection{Tessellation of the Study Area}
We partition the study area using a regular tessellation, where the geometric center of each tile is defined as a centroid. Let $\pazocal{C}$ denote the set of all centroids in the study area. A traveler can choose from various transportation modes for a trip between any pair of centroids. For instance, a traveler may walk the entire distance at speed $v_{walking}$, or they could walk from the origin centroid $c_i$ to a PT station, then take PT to another station, and finally walk from there to the destination centroid $c_j$. Observe that while actual trips may be longer or shorter than centroid-to-centroid distances, these variations tend to balance out across the PT network. This is why we approximate all travel as centroid-to-centroid.

\subsection{Graph Model of PT Network}

We model the original PT structure before disruption as graph $\pazocal{G}^o = (\pazocal{V}^o, \pazocal{E}^o)$. $\pazocal{G}^o$ is composed of multiple PT lines. Each line $\ell$ has a headway $t_\ell$, defined as the distance between two consecutive vehicles expressed in time~\cite{vuchic1981timed}. The headway can be calculated using the formula $t_\ell = \frac{T_\ell}{N_\ell}$, where $T_\ell$ represents the round trip time along the entire line, and $N_\ell$ denotes the number of vehicles operating on that line~(Eq.~4.4 of \cite{vuchic1981timed}). A line consists of a sequence of stations connected by edges. 
The average waiting time for line $\ell$ is $\frac{t_\ell}{2}$~\cite{pinto2020joint}, We define $t(s_j, s_{j+1})$, the time taken to travel between any two successive stations $s_j$ and $s_{j+1}$, and $t_{s_j}$, the dwell time at station $s_j$. For interchange at station $s_j$ from line $\ell$ to line $\ell'$, we consider an average waiting time $\frac{t_{\ell'}}{2}$ at station $s_j$.

The set of centroids and edges connecting any centroid to every other centroid and every station is likewise included in $\pazocal{G}^o$. Thus, the set of nodes $\pazocal{V}^o$ is:
$\pazocal{V}^o = \pazocal{S} \cup \pazocal{D} \cup \pazocal{C}$
where $\pazocal{S}$ represents bus stops, $\pazocal{D}$ urban rail transit stations, and $\pazocal{C}$ centroids in the study area.
The set of edges $\pazocal{E}^o$ represents the connections between these nodes. 

We consider in this case the disruption of an entire rail line, which often occurs for maintenance work~\cite{christoforou2016managing,20minutes2024}. In this case, affected rail stations are deactivated, and the edges connecting them are unusable. An example of this situation is depicted in Fig.~\ref{fig:PT_repre}, where the gray thick line represents the disrupted urban rail transit (URT) line and stations of this line, indicated by circles with red crosses at the center, are closed (We will come back to the other illustration within this figure in section ``A two-stage heuristic algorithm''). From the modeling point of view, the PT network with the disrupted line is represented by a new PT graph, $\pazocal{G}^\text{DISR}$. Let $\pazocal{D}$ denote the set of disrupted stations; we assume that a central system advises passengers affected by disruptions in stations $d\in\pazocal{D}$ to board redesigned buses at particular consolidation nodes (by phone messages or signage). The consolidation node for a disrupted rail station $d \in \pazocal{D}$ can be the bus stop in $\pazocal{S}$ or disrupted rail station $d$. Let $\pazocal{N}$ be the set of consolidation nodes. 

\subsection{Problem Statement}

In the face of urban rail system disruptions, our research aims to maximize accessibility through a redesign of bus lines. We extend operating bus routes allowing buses of one line to serve both stops of the original line and consolidation points at which passengers concerned by the disruption can be picked up and dropped off. 
Our optimization problem makes the following five decisions:
(a) To which candidate consolidation points should passengers from the disrupted rail stations be suggested to go?
(b) To which bus lines should the consolidation nodes be assigned?
(c) Which bus line should be extended in order to serve such passengers?
(d) When a bus serves multiple consolidation points, what is the optimal order of visits?
(e) How to re-assign the existing fleet of buses and allocate additional vehicles (if available) among the original and extended lines?

The practical implementation of this approach assumes that drivers can adapt to revised routes through real-time navigation systems commonly employed by transit operators worldwide. These systems, such as onboard consoles or tablets, provide turn-by-turn navigation and can rapidly deploy route changes during disruptions with minimal additional hardware requirements.

\begin{figure}[!t]
    \centering
    \includegraphics[width=1\linewidth]{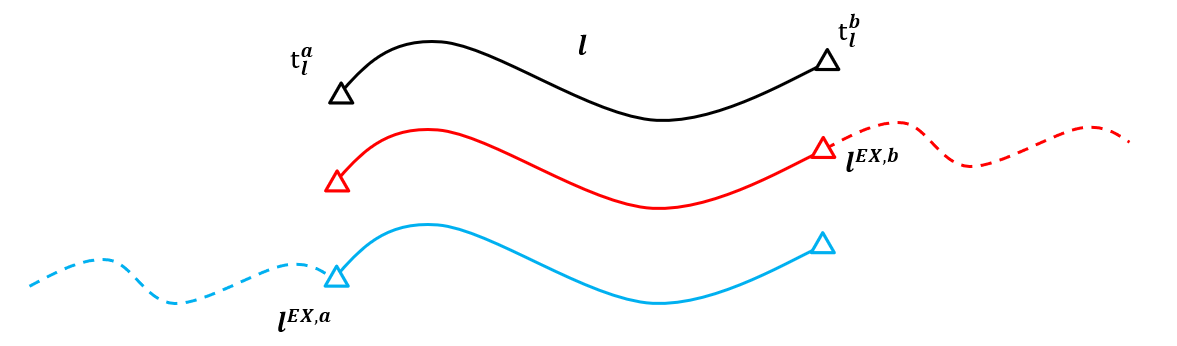}
    \caption{Line extension illustration}
    \label{fig:lineEX}
\end{figure}

As illustrated in Fig.~\ref{fig:lineEX}, for any given line $\ell$, let $t_{\ell}^a$ and $t_{\ell}^b$ represent its terminals. We can activate an extended version of $\ell$ from terminal $t_{\ell}^b$, which will serve the stops of the regular $\ell$ as well as consolidation nodes that will be associated with this line. Let’s denote this extended line as $\ell^{EX, b}$ (red line). Similarly, we can define $\ell^{EX, a}$ (blue line). Let $\pazocal{L}_\text{REG} = \{\ell_1, \ell_2, ..., \ell_\ell\}$ be the set of regular lines. Let $\pazocal{L}_\text{EX} = \{\ell^{EX, a}, \ell^{EX, b} | \ell \in \pazocal{L}_\text{REG}\}$ the set of all potential extended bus lines.

\begin{table}[t]
\centering
\caption{Mathematical notation for the optimization model.}
\label{tab:model-notations}
\footnotesize
\begin{tabularx}{\columnwidth}{@{}l>{\raggedright\arraybackslash}X@{}}
\toprule
Symbol & Description \\
\midrule
$n$ & Consolidation node index. \\
$d$ & Demand point index. \\
$i, j$ & Indices for nodes. \\
$\ell$ & Bus line index. \\
\midrule
Sets & $\pazocal{V} = \pazocal{N} \cup \pazocal{D} \cup \pazocal{C}$. \\
\midrule
$\pazocal{N}$ & Set of consolidation nodes, i.e., boarding locations for affected passengers. \\
$\pazocal{D}$ & Set of disrupted rail stations, where affected passengers are gathered. \\
$\pazocal{C}$ & Set of centroids in the study area. \\
$\pazocal{L} = \pazocal{L}_\text{REG} \cup \pazocal{L}_\text{EX}$ & Set of all bus lines in the network, where $\pazocal{L}_\text{REG}$ are regular lines and $\pazocal{L}_\text{EX}$ are extended lines proposed by the model. \\
\midrule
Parameters & \\
\midrule
$d_{ij}$ & Distance between nodes $i$ and $j$, $\forall i,j \in \pazocal{V}$. \\
$q_d$ & Passenger demand at disrupted rail station $d$, $\forall d \in \pazocal{D}$. \\
$cap$ & Passenger capacity per bus. \\
$N^o_\ell$ & Number of buses on regular line $\ell$, $\forall \ell \in \pazocal{L}_\text{REG}$. \\
$N_{max}$ & Maximum number of vehicles that can be added. \\
$D_{max}$ & Maximum distance allowed from disrupted station to consolidation node. \\
\midrule
Decision Variables & \\
\midrule
$\pazocal{G}$ & PT structure. \\
$n_\ell \in \mathbb{Z}^+$ & Additional vehicles per unit time on line $\ell$ or its extensions, $\forall \ell \in \pazocal{L}_\text{REG}$. \\
$x_\ell \in \mathbb{Z}^+$ & Vehicles per unit time on line $\ell$, $\forall \ell \in \pazocal{L}$. \\
$y_{\ell^{EX}} \in \{0,1\}$ & 1 if extended line $\ell^{EX} \in \pazocal{L}_\text{EX}$ is activated, 0 otherwise. \\
$w_{n\ell} \in \{0,1\}$ & 1 if consolidation point $n$ is assigned to extended line $\ell$, 0 otherwise. \\
$z_{dn} \in \{0,1\}$ & 1 if disrupted station $d$ is assigned to consolidation $n$, 0 otherwise. \\
$u_i^{\ell} \in \mathbb{Z}^+$ & Positional order of node $i$ in line $\ell$ (0 if not visited). \\ 
$X_{ij}^{\ell} \in \{0,1\}$ & 1 if line $\ell$ connects $i$ and $j$, 0 otherwise. \\
\bottomrule
\end{tabularx}
\end{table}

\subsection{Problem Formulation}
\subsubsection{Objective function}
Following the gravity-based definition of accessibility~\cite{miller2020measuring}, we interpret accessibility as the number of opportunities reachable within a given travel time. We make the assumption that all trips start and end at the centroids. Accessibility of a centroid $c_i \in \pazocal{C}$ is then given by:
\begin{equation}
\text{acc}(c_i, \pazocal{G}) = \sum_{c_j \in \pazocal{C}, c_j \neq c_i} \frac{O_{c_j}}{T_{\pazocal{G}}(c_i, c_j)}\label{equ::accINIT}
\end{equation}
where 
$O_{c_j}$ is the number of opportunities within the tile with centroid $c_j$ and $T_{\pazocal{G}}(c_i, c_j)$ is the shortest travel time between centroids $c_i$ and $c_j$, considering that trips can combine multiple lines and can also include walk from $c_i$ to $c_j$ or to/from some intermediate stops; Subscript $\pazocal{G}$ indicates that travel time depends on PT structure $\pazocal{G}$. 
To implement this accessibility measure, we partition the study area using a regular tessellation with 1-kilometer grid cells, where centroids are defined as geometric centers of each tile. The 1-kilometer grid size balances computational efficiency with spatial precision, representing a reasonable walking distance while capturing accessibility variations in suburban areas with sparse public transport networks. Note that $O_{c_j}$ is an input parameter, exogenously determined.
Indeed, if the PT structure changes (for instance some edges are no more available because of disruption), some shortest paths may change accordingly, resulting in a change in travel time~$T_\pazocal G(c_i, c_j)$.
Let $\pazocal{G}^o$ be the original PT structure before disruption, and $\pazocal{G}^\text{DISR}$ be the disrupted PT structure. Then, it's easy to show that:
\begin{align}
&& 
\text{acc}(c_i, \pazocal{G}^\text{DISR}) \leq \text{acc}(c_i, \pazocal{G}^{o})&& \forall c_i \in \pazocal{C}
\end{align}

Our problem consists in finding a new PT structure $\pazocal{G}$ such as to maximize the overall accessibility index:
\begin{align}
&& 
\max_{\pazocal{G}}\quad\sum_{c_i \in \pazocal{C}}\text{acc}(c_i, \pazocal{G}) = \max_{\pazocal{G}}\quad\sum_{c_i \in \pazocal{C}}\sum_{c_j \in \pazocal{C}, c_j \neq c_i} \frac{O_{c_j}}{T_{\pazocal{G}}(c_i, c_j)}
\label{equ::acc}
\end{align}

\subsubsection{Constraints}

\begin{align}
& \text{Demand allocation constraints:}\nonumber \\
& \sum_{\ell \in \pazocal{L}_\text{EX}} w_{n\ell} \leq 1, && \forall n \in \pazocal{N} \label{eq:1} \\
& \sum_{n \in \pazocal{N}} z_{dn} = 1, && \forall d \in \pazocal{D} \label{eq:2} \\
& z_{dn} \leq \sum_{\ell \in \pazocal{L}_\text{EX}} w_{n\ell}, && \forall n \in \pazocal{N}, \forall d \in \pazocal{D} \label{eq:3} \\
& z_{dn} \cdot d_{nd} \leq D_{max}, && \forall n \in \pazocal{N}, \forall d \in \pazocal{D} \label{eq:4}
\end{align}
Constraints Eq.~\ref{eq:1} to Eq.~\ref{eq:4} handle the allocation relationship between disrupted rail stations and consolidation nodes: Constraint Eq.~\ref{eq:1} ensures each consolidation node is assigned to at most one extended line; 
Constraint Eq.~\ref{eq:2} guarantees each disrupted rail station is allocated to exactly one consolidation node; 
Constraint Eq.~\ref{eq:3} ensures disrupted rail stations are only assigned to stops served by extended lines; 
Constraint Eq.~\ref{eq:4} limits the distance between a disrupted rail station and its assigned consolidation node to within the maximum allowable distance.

\begin{align}
& \text{Line extension routing constraints:}\nonumber \\
&\begin{aligned}
   &\sum_{j \in \pazocal{N}} X_{t_\ell^aj}^{\ell^{EX,a}} = y_{\ell^{EX,a}};\\ &\sum_{j \in \pazocal{N}} X_{t_\ell^bj}^{\ell^{EX,b}} = y_{\ell^{EX,b}} 
\end{aligned} && \forall \ell \in \pazocal{L}_\text{REG}\label{Seq:1} \\
&\sum_{j \in \pazocal{N}} X_{ij}^{\ell} = \sum_{j \in \pazocal{N}} X_{ji}^l, && \forall i \in \pazocal{N}, \forall \ell \in \pazocal{L}_\text{EX}\label{Seq:3} \\
&\sum_{j \in \pazocal{N}} X_{ij}^{\ell} \leq w_{i\ell}, && \forall i \in \pazocal{N}, \forall \ell \in \pazocal{L}_\text{EX}\label{Seq:4} \\
&\sum_{i \in \pazocal{N}} X_{ij}^{\ell} \leq w_{j\ell}, && \forall j \in \pazocal{N}, \forall \ell \in \pazocal{L}_\text{EX}\label{Seq:5} \\
&\begin{aligned}
    &u_i^{\ell} - u_j^\ell + (|\pazocal{N}|+2) \\
    &\quad\cdot X_{ij}^{\ell} \leq (|\pazocal{N}|+1) \cdot w_{i\ell}
\end{aligned},&&  \begin{aligned}
    &\forall i,j \in \pazocal{N}, i \neq j, \\
    &\forall \ell \in \pazocal{L}_\text{EX}
\end{aligned}\label{Seq:6} \\
&2 \cdot w_{i\ell} \leq u_i^{\ell} \leq (|\pazocal{N}|+1) \cdot w_{i\ell},&&  \forall i \in \pazocal{N}, \forall \ell \in \pazocal{L}_\text{EX}\label{Seq:7} \\
&u_{t_\ell^a}^{\ell^{EX,a}} = y_{\ell^{EX,a}}~; u_{t_\ell^b}^{\ell^{EX,b}} = y_{\ell^{EX,b}}&& \forall \ell \in \pazocal{L}_\text{REG}\label{Seq:8}
\end{align}

Constraints Eq.~\ref{Seq:1} to Eq.~\ref{Seq:8} design the line extensions: 
Constraint Eq.~\ref{Seq:1} defines the starting point constraints for line extensions. As shown in Figure~\ref{fig:lineEX}, each regular line $\ell$ has two terminals $t_\ell^a$ and $t_\ell^b$. The first term $\sum_{j \in \pazocal{N}} X_{t_\ell^aj}^{\ell^{EX,a}} = y_{\ell^{EX,a}}$ ensures that if line extension from terminal $t_\ell^a$ is activated, there must exist exactly one direct connection from $t_\ell^a$ to a consolidation node. Similarly, the second term $\sum_{j \in \pazocal{N}} X_{t_\ell^bj}^{\ell^{EX,b}} = y_{\ell^{EX,b}}$ addresses the extension from terminal $t_\ell^b$.
Constraint Eq.~\ref{Seq:3} maintains flow conservation at consolidation nodes and ensures that the number of incoming and outgoing connections is equal. 
Constraints Eq.~\ref{Seq:4} and Eq.~\ref{Seq:5} limit each node to have at most one outgoing and one incoming edge in each line, respectively, preserving the linear structure of the line. 
Constraint Eq.~\ref{Seq:6} implements the Miller-Tucker-Zemlin subtour elimination constraint~\cite{miller1960integer}, ensuring that the line forms a single path without cycles among consolidation nodes. Lastly, constraint Eq.~\ref{Seq:7} sets the bounds for node position variables, maintaining sequential order, while Constraint Eq.~\ref{Seq:8} specifically designates the starting terminal. 

\begin{align}
& \text{Vehicle allocation and line selection constraints: }\nonumber \\
& x_{\ell^{EX,a}} \geq y_{\ell^{EX,a}}~; x_{\ell^{EX,b}} \geq y_{\ell^{EX,b}}&& \forall \ell \in \pazocal{L}_\text{REG} \label{eq:9} \\
& x_{\ell^{EX,a}} \leq M \cdot y_{\ell^{EX,a}}~; x_{\ell^{EX,b}} \leq M \cdot y_{\ell^{EX,b}}&& \forall \ell \in \pazocal{L}_\text{REG} \label{eq:10} \\
& cap \cdot x_{\ell} \geq \sum_{d \in \pazocal{D}} \sum_{n \in \pazocal{N}} q_d \cdot z_{dn} \cdot w_{n\ell}, && \forall \ell \in \pazocal{L}_\text{EX} \label{eq:11} \\
& n_\ell = x_\ell + x_{\ell^{EX,a}} + x_{\ell^{EX,b}} - N^o_\ell, && \forall \ell \in \pazocal{L}_\text{REG}\label{eq:12} \\
& \sum_{\ell \in \pazocal{L}} n_\ell \leq N_{max} \label{eq:14}
\end{align}
Constraints Eq.~\ref{eq:9} to Eq.~\ref{eq:14} deal with vehicle allocation and line selection: 
Constraint Eq.~\ref{eq:9} and Constraint Eq.~\ref{eq:10} ensure vehicles are only assigned to a line extension if this extended line is activated; 
Constraint Eq.~\ref{eq:11} guarantees the vehicle capacity of each extended line meets demand; Note that $x_{\ell}$ indicates the number of vehicles operating on the extended line $\ell$ and $z_{dn} \cdot w_{n\ell}$ is one if and only if disrupted rail station $d$ is served by extended line $\ell$;
$q_d$ can be estimated through analysis of historical passenger flow data and real-time monitoring etc. In our numerical experiments, we employed statistical models; this constraint makes the problem non-linear.
Constraint Eq.~\ref{eq:12} governs the allocation of vehicles across the transit network. For each regular line $\ell$, it ensures that the number of vehicles on the regular line ($x_\ell$) equals the initial fleet size ($N^o_\ell$) plus additional vehicles fleet size ($n_\ell$), minus the number of vehicles on extended lines. 
Constraint Eq.~\ref{eq:14} expresses a budget for additional buses to be added. 
Note that the formulation allows us to reduce the number of buses operating on regular lines in our solution with respect to the PT design before disruption. It may indeed be useful to slightly reduce the service on regular lines to compensate for the loss of accessibility due to disruption.

\subsection{Considerations about Computational Complexity}
Our problem can be viewed as a complex combination of the Vehicle Routing Problem (VRP)~\cite{toth2002vehicle} and the Capacitated Facility Location Problem (CFLP)~\cite{wu2006capacitated}, with an even more sophisticated objective function. Constraints Eq.~\ref{eq:9}~- Eq.~\ref{eq:14} resemble those in VRP, handling vehicle allocation and line selection, while constraints Eq.~\ref{eq:1}~- Eq.~\ref{eq:4} and parts of constraints Eq.~\ref{Seq:1}~- Eq.~\ref{Seq:8} mirror CFLP constraints, addressing consolidation node selection, disrupted rail station assignment, and capacity limitations. Our objective function (Eq.~\ref{equ::acc}) is more complex than those in VRP or CFLP; additionally, constraints  Eq.~\ref{Seq:1}~- Eq.~\ref{Seq:8} incorporate line design elements, further increasing the problem's complexity. Given that both VRP and CFLP are known NP-hard problems, and our problem not only combines these two sub-problems but also adds extra constraints and a more intricate objective function, we can reasonably conclude that our problem is NP-hard.

\section{A Two-stage Heuristic Algorithm}
\begin{algorithm}
\caption{Two-Stage resolution method}
\label{alg:twoSTAGE}
\begin{algorithmic}[1]
\STATE \textbf{Stage 1: Demand-Driven Consolidation and Route Optimization}
\STATE Assign disrupted rail station $d\in\pazocal{D}$ to candidate consolidation nodes $n\in\pazocal{N}$ according to Eq.~\ref{eq::defZ}, cluster above-assigned consolidation nodes $n$ into clusters $\pazocal{K}$.
\FOR{each cluster $k \in \pazocal{K}$}
    \FOR{each regular line $\ell \in \pazocal{L}_\text{REG}$}
    \STATE (Refer to Fig.~\ref{fig:lineEX}) Compute extended line {$\ell_k^{EX, a}$,} originating from terminal $t^a_\ell$ of line $\ell$ and traversing all consolidation nodes $n$ within cluster $k${. To compute line $\ell_k^{EX, a}$, just apply the shortest path among nodes~$n$ in cluster~$k$.} {Perform the same computation to find extended line} line {$\ell_k^{EX, b}$}.
  \ENDFOR
\ENDFOR
\STATE \textbf{Stage 2: Bus Routes and Vehicle Allocation Optimization}
\STATE Decide which line should be connected to which cluster via line extensions the number of buses in each line (regular and extended) via the procedure explained in subsection ``Stage~2''.
\RETURN Optimized redesigned PT network
\end{algorithmic}
\end{algorithm}

{To address the computational complexity of our NP-hard optimization problem, we propose a two-stage heuristic algorithm that decomposes the computational burden while preserving problem structure integrity.}
In stage 1, we employ spatial clustering to group consolidation nodes into a reduced set of clusters and calculate optimal routing within each cluster. In stage 2, each line is associated to a cluster by one of its bus terminals, then extends through this terminal and all consolidation nodes of the associated cluster, following the routing decided in stage 1.

{As shown in Fig.~\ref{fig:PT_repre}, our algorithm addresses a scenario where a URT line (depicted by the thick grey line) is disrupted, causing six stations (red dots) to cease operations, necessitating a redesign of the bus service to restore the loss of accessibility. The process begins by guiding affected passengers to the suggested consolidation nodes, which are defined in Eq.~\ref{eq::defZ}; in the example of Fig.~\ref{fig:PT_repre}, the clusters output from stage 1 are represented by red dashed circles.} Let us focus on the cluster composed of consolidation nodes $n_1$ and $n_2$; suppose stage 2 determines that the line terminating at this terminal $t_1$ is associated with the above-mentioned cluster. A new extension line is thus created from $t_1$, shown as the dotted line; this means that bus line 4 now serves the original bus stops plus consolidation nodes $n_1$ and $n_2$. Tab.~\ref{tab::model-stage2} summarizes the clustering-related parameters and variables.

\begin{table}
\centering
\caption{Notation for Cluster-Based Formulation}
\begin{small}
\begin{tabular}{l>{\raggedright\arraybackslash}p{5cm}}
\toprule
Symbol\&Sets & Description\\
\midrule
$k$ & Cluster index, $\forall k \in \pazocal{K}$.\\
$\pazocal{K}$ & The set of clusters. \\
\midrule
Parameters & \\
\midrule
$d_{\ell_k^{\text{EX},a}}~; d_{\ell_k^{\text{EX},b}}$ & The shortest distance of the hypothetical extended lines traversing the consolidation nodes within cluster $k$, $\forall \ell \in \pazocal{L}_\text{REG}, \forall k \in \pazocal{K} $. \\
$q_k$ & Number of passenger demands per unit of time within cluster $k$, $\forall k \in \pazocal{K}$. \\
$acc_\ell$ & Proxy measure of line $\ell$'s importance for accessibility, $\forall \ell \in \pazocal{L}_\text{REG}$ \\
$acc_n$ & Proxy measure of node $n$'s importance for accessibility, $\forall n \in \pazocal{N}$ \\
$acc^k$ & Proxy measure of cluster $k$'s importance for accessibility, $\forall k \in \pazocal{K}$ \\
$acc_\ell^k$ & Proxy measure of extended line $\ell$'s importance for accessibility when serving cluster $k$, $\forall \ell \in \pazocal{L}_\text{REG}, \forall k \in \pazocal{K}$ \\
\midrule
Decision Variables & \\
\midrule
$y_{\ell k} \in {0, 1}$ & 1 if extended line $\ell$ serves cluster k, 0 otherwise, $\forall \ell \in \pazocal{L}_\text{EX}, \forall k \in \pazocal{K}$. \\
\bottomrule
\end{tabular}\label{tab::model-stage2}
\end{small}
\end{table}
\subsection{Stage 1: Demand-Driven Consolidation and Route Optimization}

The first stage involves 
(i) associating disrupted rail stations to consolidation nodes based on constraint Eq.~\ref{eq:1}~- Eq.~\ref{eq:4},
(ii) clustering consolidation nodes, 
(iii) optimizing the traversal order for the extended lines of the suggested consolidation nodes of cluster $k$ to serve as an extended line $\ell$. 
To simplify, we associate rail station $d$ with consolidation node $n_d$, chosen as follows: if there exist consolidation nodes in $\pazocal{S}$ within distance $D_{max}$ from $d$, we set $n_d$ as the closest among them; otherwise, $n_d = d$, i.e. the consolidation node associated to the disrupted rail station is the station itself. 
We preferentially choose existing bus stops in $\pazocal{S}$ to consolidate disrupted demand in order to exploit the facilities available there (e.g., shelters, dynamic indication screens). Moreover, passengers can also board regular lines once they arrive at an existing bus stop.
Therefore,
\begin{equation}
z_{nd} = 
\begin{cases} 
1 & \text{if } n = n_d \\ 
0 & \text{otherwise}
\end{cases}\label{eq::defZ}
\end{equation}
In order for neighboring consolidation nodes to be served by the same bus routes and minimize the number of extended bus routes to avoid more impact on the regular line, we employ DBSCAN for clustering due to its ability to autonomously determine the number of clusters~\cite{ester1996density}. For each resulting cluster $k \in \pazocal{K}$ and for any regular line~$\ell \in \pazocal L$, we compute a candidate extension line~$\ell_k^{\text{EX}, a}$, which represents the exact path that would be followed by extended line~$\ell^{\text{EX},a}$ if it were associated to cluster~$k$. To calculate this path, we formulate an open path problem based on constraints Eq.~\ref{Seq:1} -  Eq.~\ref{Seq:8}, calculating the minimum distance to traverse the consolidation nodes $n \in \pazocal{N}$ within cluster $k$. We do the same for extended line~$\ell^{\text{EX},b}$. We repeat this calculation for each extension of all regular lines. Which extended line will be actually activated and traversing which cluster will be decided in Stage~2. We denote the length of hypothetical lines~$\ell_k^{\text{EX},a},\ell_k^{\text{EX},b}$ by~$d_{\ell_k^{\text{EX},a}}, d_{\ell_k^{\text{EX},b}}$, respectively.

\subsection{Stage 2: Bus Route Extension and Vehicle Allocation Optimization}
In stage 2, we decide which extended line to associate to which cluster via line extension (variable~$y_{\ell k}, \forall k \in \pazocal{K}$) and show how many buses per unit of time we should allocate to each extended line and regular line (variable~$x_\ell$, $\forall \ell \in \pazocal{L} = \pazocal{L}_{\text{REG}} \cup \pazocal{L}_{\text{EX}}$). 
Such decisions should be taken in order to maximize the objective~(Eq.~\ref{equ::acc}).

However, directly incorporating the accessibility calculation (Eq.~\ref{equ::acc}) into the IP model presents significant challenges. Indeed, accessibility calculation requires the prior computation of one-to-many shortest paths from any location; in large-scale public transport networks, calculating the shortest travel times involving multiple line combinations and walking options leads to a combinatorial explosion, significantly increasing model complexity and solution time.

For this reason, we replace the objective in Eq.~\ref{equ::acc} with a surrogate objective~$f=f_1-f_2$, where $f_1$ is a proxy of overall accessibility (Eq.~\ref{equ::acc}) and~$f_2$ is the distance covered by extension buses. 
{Our two-stage heuristic algorithm effectively reduces computational complexity where only $acc_\ell$ (proxy measure of regular bus line $\ell$'s importance for accessibility,  Eq.~\ref{eq::accl}) and $acc_n$ (proxy measure of consolidation node $n$'s importance for accessibility,  Eq.~\ref{eq::accn})) need to be pre-computed as input data.
These accessibility scores are used directly in the calculation of our surrogate objective function Eq.~\ref{eq::surF1}. Once these metrics are available, our solution can be computed within seconds, making it suitable for responding to disruptions occurring on a timescale of hours.}
To calculate~$f_1$, instead of computing the accessibility of all centroids and for each alternative, PT re-designed the configuration; we pre-computed the accessibility from stops and consolidation nodes in order to avoid combinatorial explosion. 
We associate regular bus line~$\ell$ with the following accessibility proxy score:

\begin{align}
&& 
\text{acc}_\ell = 
\frac{1}{|\pazocal{S}_\ell|} 
\sum_{s \in \pazocal{S}_\ell} 
\sum_{c_j \in \pazocal{C}, c_i \neq s} 
\frac{O_{c_j}}
{T_{\pazocal{G}^\text{DISR}}(s, c_j)}
, && \forall l\in\pazocal L_\text{REG}
\label{eq::accl}
\end{align}
where $S_\ell$ is the set of all stations on regular line $\ell$, $s$ is a stop on the regular line $\ell$, $O_{c_j}$ is the number of opportunities within the tile with centroid $c_j$, {note that $O_{c_j}$ is an input parameter, exogenously determined (i.e., we do not have control on it)''.}
$T_\pazocal{G}^\text{DISR}(s, c_j)$ is the travel time from stops $s$ to centroid $c_j$ in the PT network structure during disruption; with slight abuse of notation, $c_j \neq s$ means that $c_j$ and $s$ must not be in the same tile.
Similarly, we can get the accessibility proxy of consolidation node $n$ as $acc_n$, as shown in Eq.~\ref{eq::accn}; where again $c_j \neq n$, indicates that $c_j$ and $n$ must not be in the same tile. We associate to cluster $k$ the accessibility proxy score as shown in Eq.~\ref{eq::acck}.
\begin{align}
&& 
\text{acc}_n = 
\sum_{c_j \in \pazocal{C}, c_i \neq n} 
\frac{O_{c_j}}
{T_{\pazocal{G}^\text{DISR}}(n, c_j)}
, && \forall n\in\pazocal{N}
\label{eq::accn}
\end{align}
\begin{align}
&& 
\text{acc}_k = 
\sum_{n \in \text{Cluter}~k}
acc_n \cdot \frac{1}{\text{Size of Cluster}~k}, && \forall k\in\pazocal{K}
\label{eq::acck}
\end{align}

 We finally associate an accessibility score with the hypothesis that an extension of the line $\ell \in l_{\text{REG}}$ will serve cluster $k \in \pazocal{K}$: 
\begin{equation}
acc_k^l = \underbrace{acc_\ell}_{\text{Eq.~}\ref{eq::accl}} + \underbrace{acc_k}_{\text{Eq.~}\ref{eq::acck}}
\end{equation}
\subsubsection{$f_1$ is defined as:}
\begin{flalign}
f_1 = \sum_{\ell \in \pazocal{L}_\text{REG}} acc_\ell \cdot x_\ell 
    + \sum_{\ell \in \pazocal{L}_\text{REG}}\sum_{k \in \pazocal{K}} acc_\ell^k \cdot (x_{\ell^{EX,a}} + x_{\ell^{EX,b}})
    \label{eq::surF1}
\end{flalign}

By maximizing~$f_1$, we favor solutions where
{
\begin{itemize}
    \item More buses are allocated to lines with stops that are more important for accessibility
    \item More buses are allocated to extended lines that serve consolidation nodes that are more important for accessibility
\end{itemize}
The rationale behind this choice is that by prioritizing lines and nodes with higher accessibility importance scores, we can more effectively restore the overall accessibility of the disrupted network.
}

\subsubsection{$f_2$ is defined as:}
\begin{flalign}
f_2 = \sum_{\ell \in \pazocal{L}_\text{REG}} \sum_{k \in \pazocal{K}} (d_{\ell^{EX,a}_k} \cdot x_{\ell^{EX,a}} + d_{\ell^{EX,b}_k} \cdot x_{\ell^{EX,b}})
\end{flalign}
By minimizing~$f_2$, we tend to preferentially allocate buses to short-distance extended lines, as we assume they are more efficient. By considering~$f_1$ and~$f_2$ together, we aim to find the delicate balance between accessibility and PT operator costs. For simplicity, we assign equal weight to these two objectives:
\begin{equation}
 \max_{y_{\ell k}, n_\ell, x_\ell} f = f_1 - f_2
\end{equation}
\subsubsection{subject to:}
\begin{align}
& \text{Constraint  Eq.}~\ref{eq:12}\nonumber \\
& \text{Constraint  Eq.}~\ref{eq:14}\nonumber \\
& q_k \leq cap \cdot \sum_{l \in L_{EX}} x_{\ell} \cdot y_{\ell k}, && \forall k \in \pazocal{K} \label{constS2::1} \\
& \sum_{\ell \in \pazocal{L}_\text{EX}} y_{\ell k} = 1, && \forall k \in \pazocal{K} \label{constS2::2} \\
& x_{\ell^{EX,a}} \geq y_{\ell^{EX,a}k}~; x_{\ell^{EX,b}} \geq y_{\ell^{EX,b}k}&& \forall \ell \in \pazocal{L}_\text{REG} \label{constS2::3} \\
& x_{\ell^{EX,a}} \leq M \cdot y_{\ell^{EX,a}k}~; x_{\ell^{EX,b}} \leq M \cdot y_{\ell^{EX,b}k}&& \forall \ell \in \pazocal{L}_\text{REG} \label{constS2::4}
\end{align}
After applying clustering, adapted from constraint Eq.~\ref{eq:9}~- Eq.~\ref{eq:14},
constraint Eq.~\ref{constS2::1} ensures that the vehicle capacity for each cluster is not exceeded, constraints Eq.~\ref{constS2::2} ensures each extended line is assigned to serve exactly one cluster; 
while constraint. Eq.~\ref{constS2::3} and  Eq.~\ref{constS2::4} ensure vehicles are only assigned in the event that an extended line is created.
\section{Numerical Results}
The performance of the methodology proposed herein is tested with numerical experiments using 2 case studies in two real cities, namely Évry-Courcouronnes, France, and Choisy-le-Roi, France. The integer programming model is coded and solved in Python using the IBM ILOG CPLEX 12.8 solver, running on an Intel Core i5 PC with 4.6 GHz speed and 16.0 GB of RAM. All solutions computed in this paper took several seconds to compute. 
\begin{table*}[!ht]
\centering
\setlength{\tabcolsep}{20pt}
\caption{Scenario parameters and hyperparameters of our resolution method}\label{tab:parameters}
    \begin{centering}
		\begin{tabular}{l l l l}
    \hline
            Name & Value & Reference  \\
    \hline
    \noalign{\vskip 0.25em} 
            Fleet size & $\dfrac{\text{Cycle time}}{\text{Headway}}$ &~Eq.~4.4 of \cite{vuchic1981timed}  \\[0.5em]
            Average waiting time & $\dfrac{headway}{2}$ & ~\cite{pinto2020joint}  \\[0.25em]
    \hline
            Length of a tile & 1 km & -  \\
            Walking speed & 3.5 km/h & Google Maps  \\
            $D_{max}$ (Eq.~\ref{eq:4}) & 500 m & -  \\
            Passenger capacity per bus (cap in Eq.~\ref{eq:11}) & 120 Passengers / Bus & ~\cite{hidalgo2013methodology}  \\
            Speed of Bus, Tram\&Metro, RER & 23.5, 35, 60 km/h & ~\cite{brtdata_paris}\cite{sortiraparis_metro}  \\
            Average headway for RER, Metro, Tram, Bus & 2, 2, 4.11, 7.17 mins & ~\cite{citytransit2024average} \\
            PT station's dwell time & 0.5 - 1 min & ~\cite{kuipers2024understanding} \\
\hline
            Gamma-Poisson mixture distribution to generate~$q_d$:\\
            \;\;\;\;\;\;\;\;Initial Shape & $k_0$ = 5 & -  \\
            \;\;\;\;\;\;\;\;Initial Scale & $s_0$ = 1 & -  \\
            \;\;\;\;\;\;\;\;Adjustment factor & $\mu$ = 1 & -  \\
            \;\;\;\;\;\;\;\;Random seed & 42 & -  \\
            Clustering hyperparamters:\\
            \;\;\;\;\;\;\;\;Optimal neighbourhood size & $eps$ = 2 km & -  \\
            \;\;\;\;\;\;\;\;Optimal minimum number of samples & $min_{samples}$ = 1 & -  \\
    \hline
		\end{tabular}
	\end{centering}
\end{table*}
\begin{figure*}[!ht]
    \centering
    \begin{subfigure}{0.43\linewidth}
        \centering
        \includegraphics[width=\linewidth]{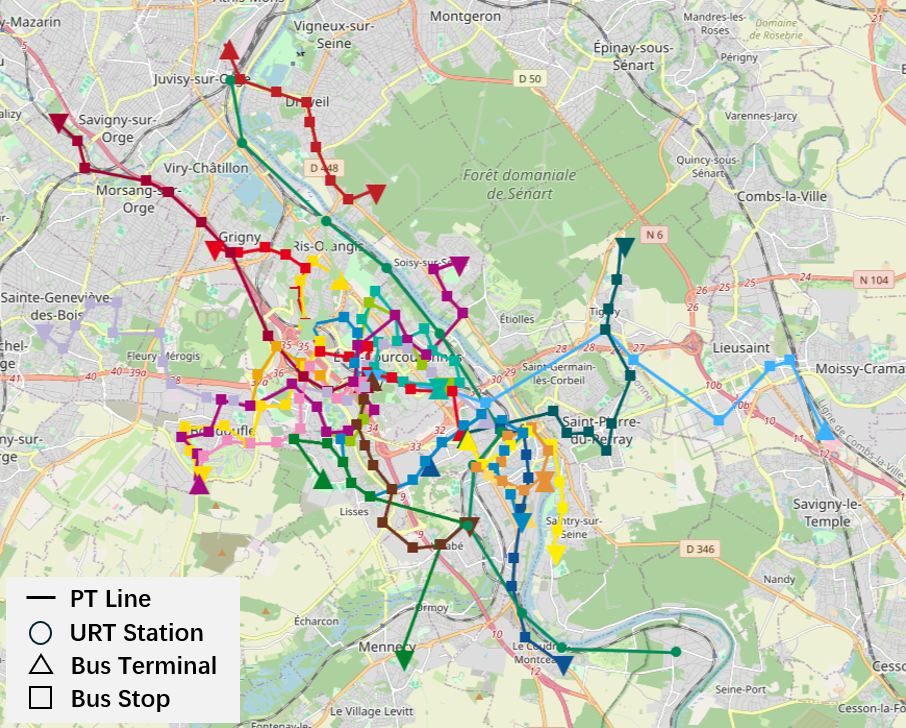}
        \caption{Map of Évry-Courcouronnes created with Open-Street Map}
        \label{fig:openstreetEVRY}
    \end{subfigure}
    \begin{subfigure}{0.46\linewidth}
        \centering
        \includegraphics[width=\linewidth]{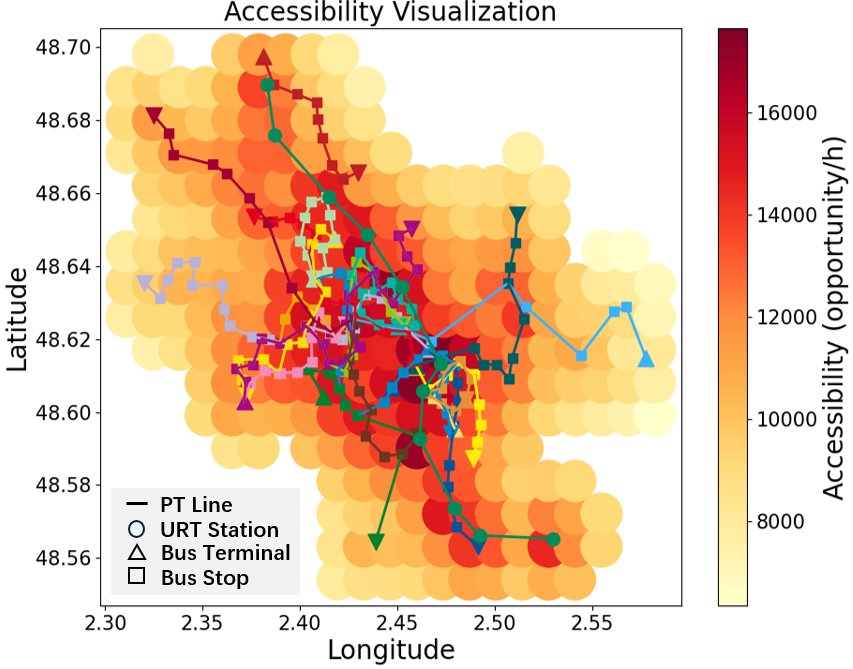}
        \caption{Heat-map of accessibility before disruption}
        \label{fig:accEVRY}
    \end{subfigure}
    \caption{{\textbf{Évry-Courcouronnes's public transport network and accessibility heat-map.} Symbols as in Fig.~\ref{fig:PT_repre}: circles = rail/tram stations, squares = bus stops, triangles = bus terminals.}}
    \label{fig:evry-transport}
\end{figure*}

\subsection{Disruption Case Study}
{We consider two study areas defined as circles with a 15 km radius centered in Évry-Courcouronnes and Choisy-le-Roi, respectively. Within these areas, we analyze public transport lines operating during daytime hours, excluding night services.}
As illustrated in Fig.~\ref{fig:openstreetEVRY}, Évry-Courcouronnes PT consists of one urban rail line (RER D), one Tramway (T12), and 19 bus lines. 
{In this figure, the different types of transport lines and stations are represented using the same symbol convention as in Fig.~\ref{fig:PT_repre}. Rail and tram stations (collectively referred to as rail transit stations) are depicted as circles, bus stops as squares, and bus terminals as triangles throughout all figures in the paper.}
Fig.~\ref{fig:accEVRY} depicts the heat map of accessibility of accessibility (Eq.~\ref{equ::acc}) where darker colors indicate higher accessibility.
A similar analysis is performed for Choisy-le-Roi. This area is served by one Tramway (T9), one urban rail line (RER C), one metro line (M8), and 16 bus lines.

In our numerical experiment, we simulated complete disruptions of the rail lines (RER D in Évry-Courcouronnes and RER C in Choisy-le-Roi) lasting several hours, resulting in the closure of all rail stations within the affected areas, dramatically reducing accessibility. To ensure accurate distance measurements for each bus line, we utilized Open-Street Map data extracted via OSMnx~\cite{boeing2017osmnx}. To compute accessibility (xwEq.~\ref{equ::accINIT}), we need first to count the number of opportunities inside each tile. The opportunities we consider are the amenities extracted from the open street map~\cite{OpenStreetMapWiki}. For simplicity, we extract all amenities without discerning among categories.

We collected the coordinates of the bus stops and urban rail stations from Île-de-France Mobilités~\cite{iledefrance2024}. 
The passenger demand is modeled using a Gamma-Poisson mixture, which calculates and averages the passenger flow for each station over seven days, resulting in different~$q_d$ values for each station~$d\in\pazocal D$. The values of headway of the considered lines are typical of peak hours.
\begin{figure*}[htbp]
    \centering
    \begin{subfigure}{0.44\linewidth}
        \centering
        \includegraphics[width=\linewidth]{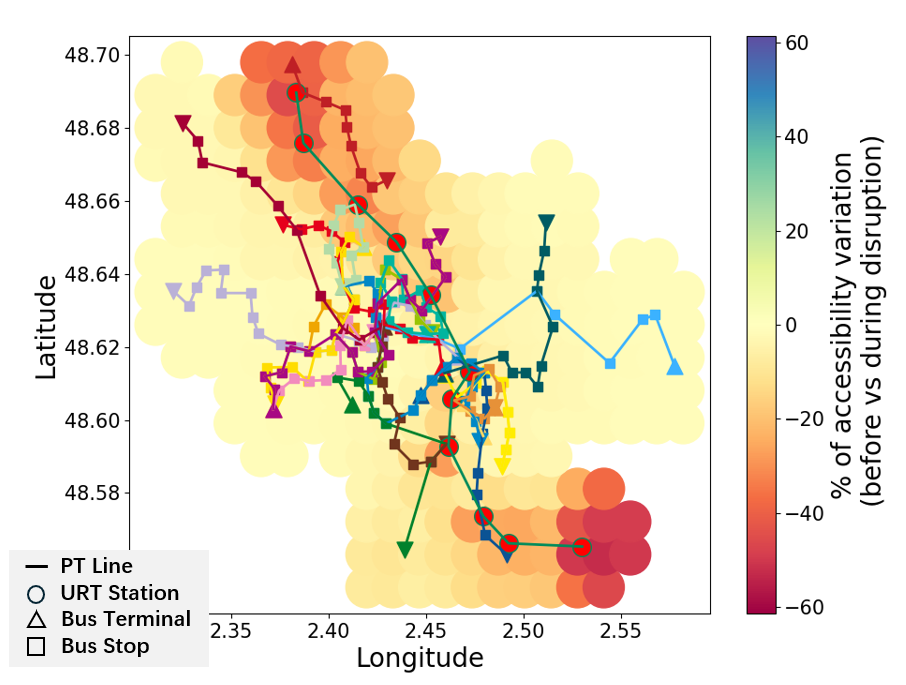}
        \caption{Impact of the disruption of the URT system (symbolic circles) in Évry-Courcouronnes region.}
        \label{fig:evryCompBDD}
    \end{subfigure}
    \begin{subfigure}{0.44\linewidth}
        \centering
        \includegraphics[width=\linewidth]{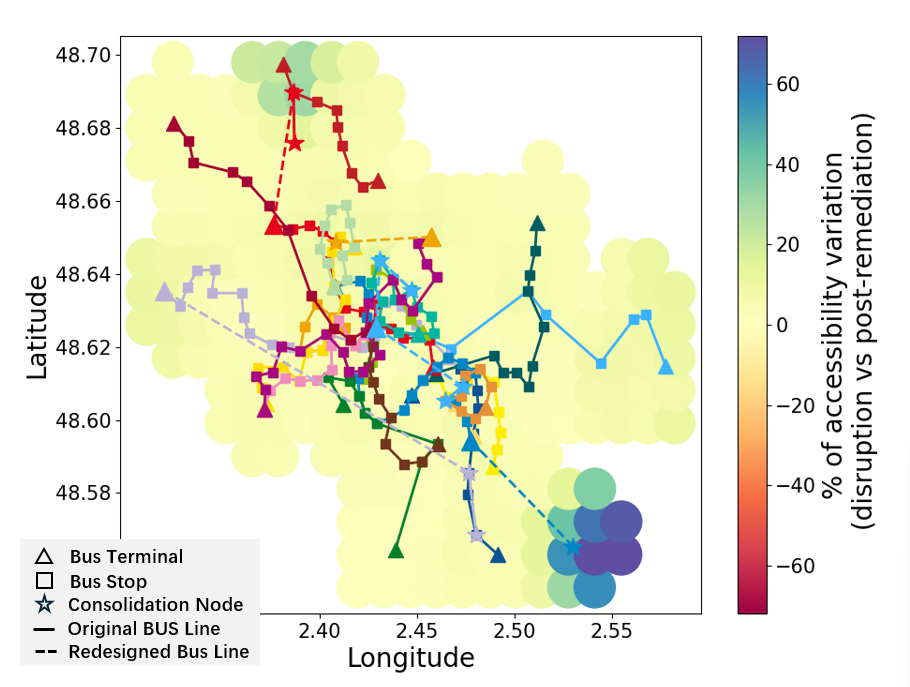}
        \caption{Accessibility recovered by our solution in Évry-Courcouronnes region (Without Adding Vehicles).}
        \label{fig:evryADD0}
    \end{subfigure}

    \begin{subfigure}{0.44\linewidth}
        \centering
        \includegraphics[width=\linewidth]{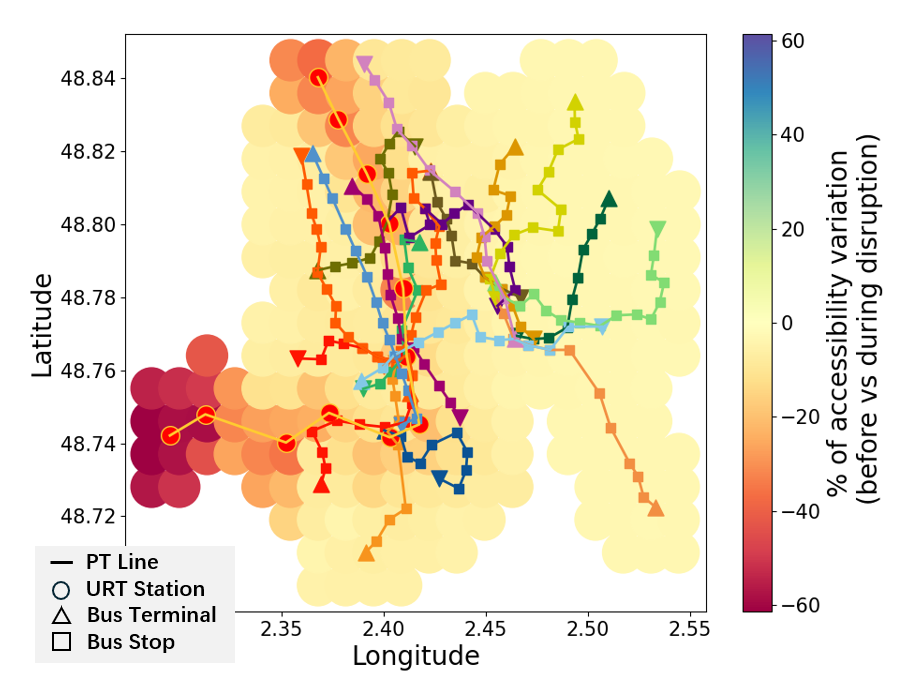}
        \caption{Impact of the disruption of the URT system (symbolic circles) in Choisy-le-Roi region.}
        \label{fig:clrCompBDD}
    \end{subfigure}
    \begin{subfigure}{0.44\linewidth}
        \centering
        \includegraphics[width=\linewidth]{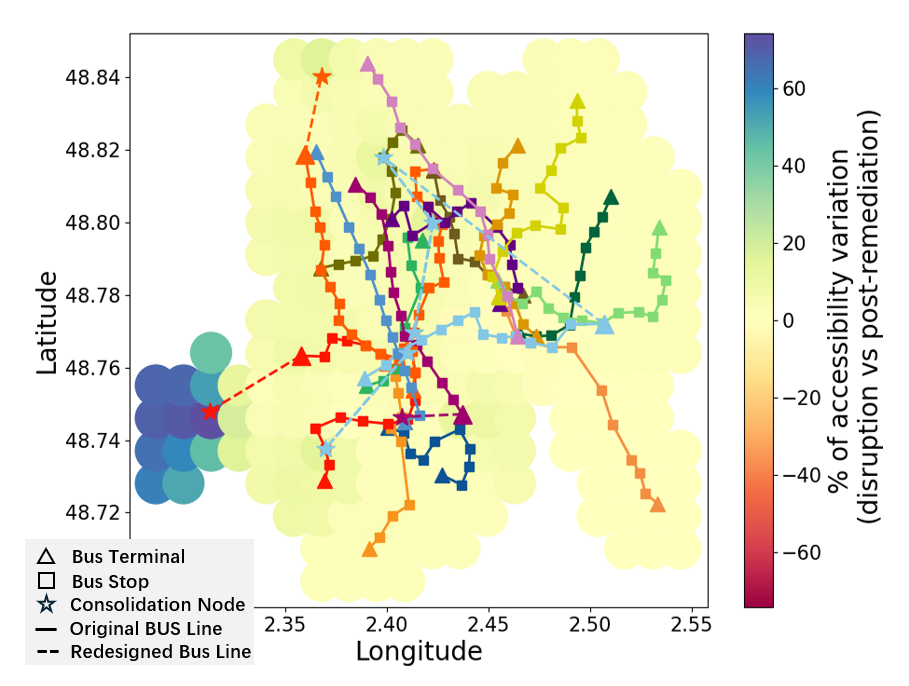}
        \caption{Accessibility recovered by our solution in Choisy-le-Roi region (Without Adding Vehicles).}
        \label{fig:clrADD0}
    \end{subfigure}
    \caption{Comparison of accessibility in two cases with no additional vehicles}
    \label{fig:accessibility_comparison_no_extra_vehicles}
\end{figure*}

All scenario parameters and algorithm hyperparameters are in Tab~\ref{tab:parameters}.

In what follows, we compare the performance of 4 scenarios:
\begin{itemize}
    \item \emph{Before Disruption}: the entire PT network is working regularly with no disruptions
    \item \emph{During Disruption}: the rail line is disrupted and no remediation has been taken yet
    \item \emph{Conventional Replacement Bus}: additional buses operate along the disrupted rail line. This is the approach classically taken. However, it is important to consider that the speed of the replacement bus is generally much lower than the replaced rail since the bus has to travel on a road network, characterized by lower speed and higher circuity
    \item \emph{Our solution}: some bus lines are extended, and the bus fleet is reallocated across regular and extended lines in order to recover the accessibility loss produced by the disruption. Such reorganization is calculated with Algorithm~\ref{alg:twoSTAGE}.
\end{itemize}

The aforementioned scenarios correspond to the following public transport (PT) network graphs: $\pazocal G^O$, $\pazocal G^\text{DISR}$, $\pazocal G^\text{REPL}$, and $\pazocal G^\text{OURS}$. The accessibility for each of these graphs are computed as acc($c_i,\pazocal G^O$), acc($c_i,\pazocal G^\text{DISR}$), acc($c_i,\pazocal G^\text{REPL}$), and acc($c_i,\pazocal G^\text{OURS}$) in accordance with Eq.~\ref{equ::accINIT}.
\subsection{Performance of Our Solution} 
Fig.~\ref{fig:evryCompBDD} and~\ref{fig:clrCompBDD}  show the variation in accessibility caused by the disruption, i.e., for each centroid $c_i$, we show
\begin{align}
&&\frac{acc(c_i,\pazocal G^\text{DISR})-acc(c_i,\pazocal G^\text{O}) }
{acc(c_i,\pazocal G^\text{O})},
&& \forall c_i\in\pazocal C    
\end{align}
Fig.~\ref{fig:evryADD0} and~\ref{fig:clrADD0}, show instead the accessibility recovered thanks to our solution, without adding any bus, i.e.:
\begin{align}
&&\frac{acc(c_i,\pazocal G^\text{OURS})-acc(c_i,\pazocal G^\text{DISR}) }
{acc(c_i,\pazocal G^\text{DISR})},
&& \forall c_i\in\pazocal C    
\end{align}
\begin{figure*}[!t]
    \centering
    \begin{subfigure}{0.47\linewidth}
        \centering
        \includegraphics[width=\linewidth]{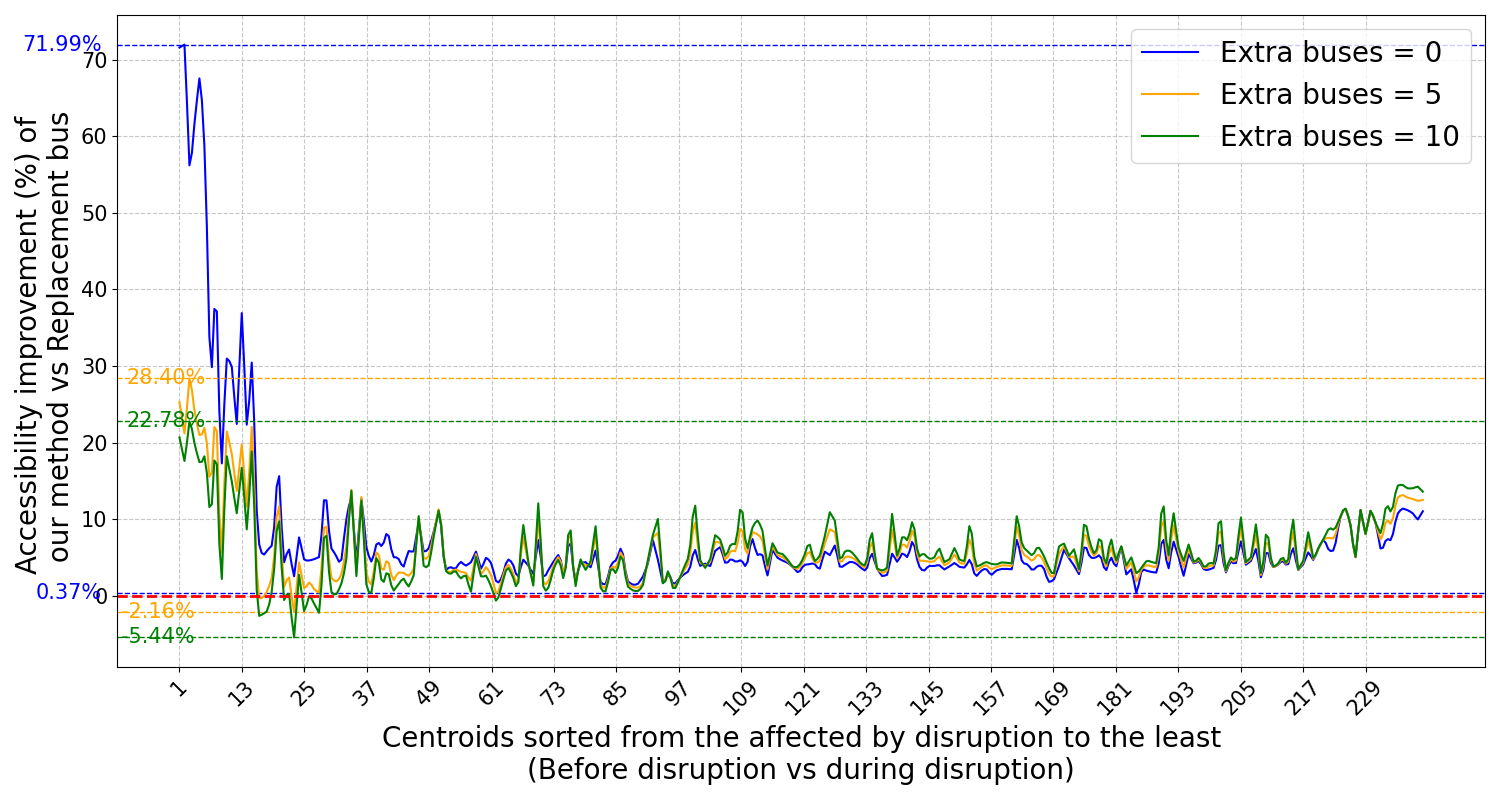}
        \caption{Évry-Courcouronnes: Impact of Extra Buses on Affected Areas.}
        \label{fig:evryACCAffected}
    \end{subfigure}
    \hfill
    \begin{subfigure}{0.47\linewidth}
        \centering
        \includegraphics[width=\linewidth]{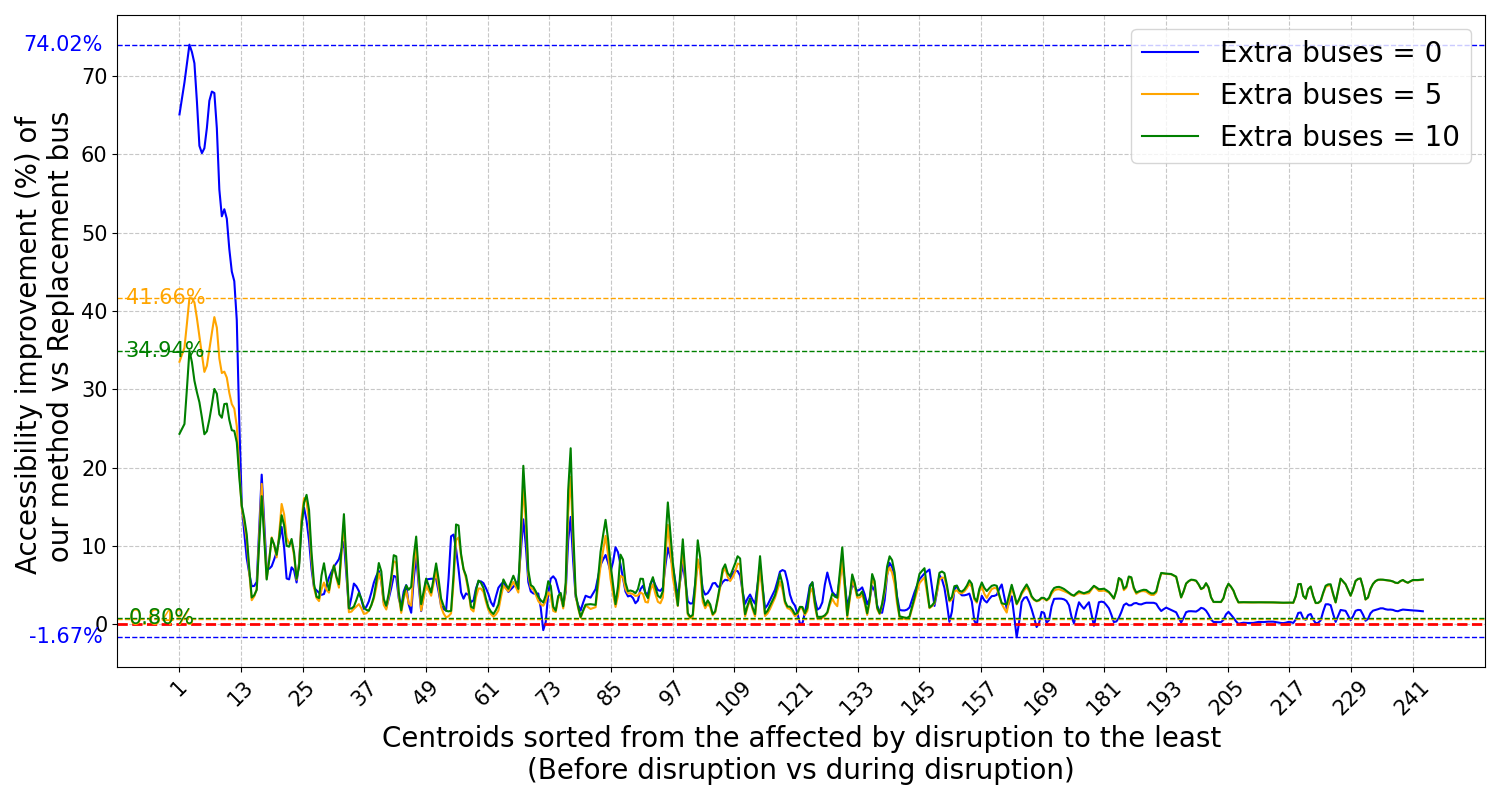}
        \caption{Choisy-le-Roi: Impact of Extra Buses on Affected Areas.}
        \label{fig:clrACCAffected}
    \end{subfigure}
    \caption{Improvement of our solution over conventional replacement bus strategy.}
    \label{fig:accCompCombined}
\end{figure*}
\begin{figure*}[!t]
    \centering
    \begin{subfigure}{0.47\linewidth}
        \centering
        \includegraphics[width=\linewidth]{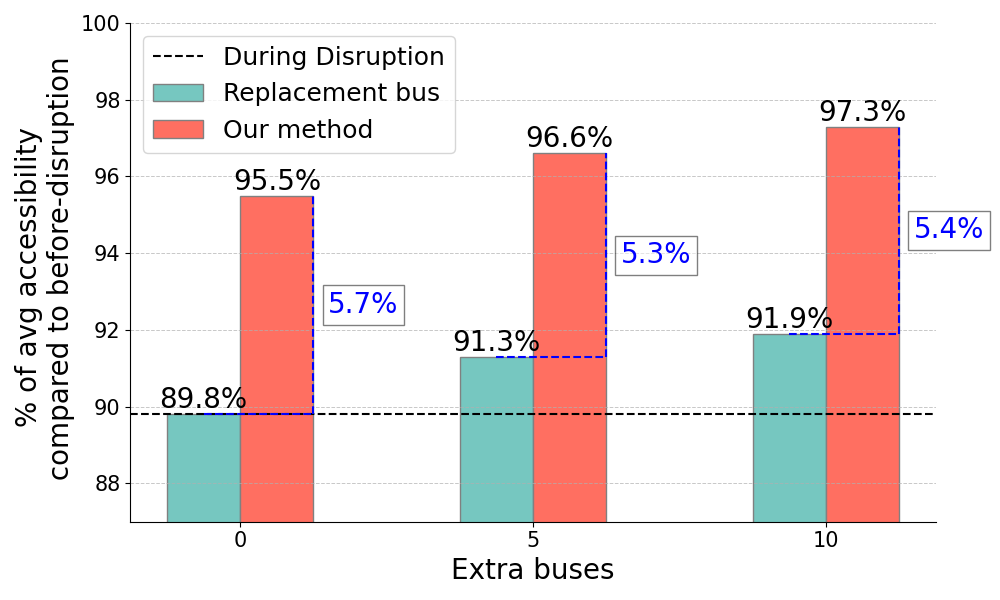}
        \caption{Évry: Average of overall accessibility}
        \label{fig:avgACCevry}
    \end{subfigure}
    \hfill
    \begin{subfigure}{0.47\linewidth}
        \centering
        \includegraphics[width=\linewidth]{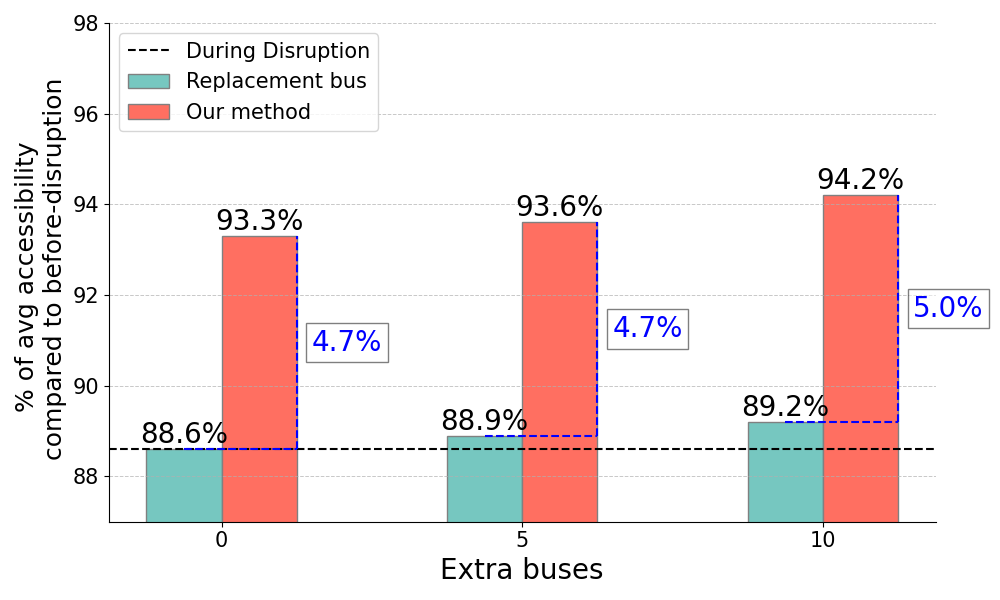}
        \caption{Choisy-le-Roi: Average of overall accessibility}
        \label{fig:avgACCclr}
    \end{subfigure}
    \caption{Accessibility recovered via remediation (values into the squares indicated the improvement of our solution over the conventional replacement bus strategy).}
    \label{fig::avgACC}
\end{figure*}

In these charts, colors closer to blue indicate an increase in accessibility, while colors closer to red signify a decrease; the symbols are the same as in Fig.~\ref{fig:PT_repre}. In Fig.~\ref{fig:evryCompBDD}, red circles are the disrupted rail stations, and in Fig.~\ref{fig:evryADD0}, the dotted lines and stars represent the routes of our re-designed bus lines. Observe that the shortest routes for the extended bus lines are computed on the road network and may not correspond to the shortest paths as the crow flies. Comparing Figures Fig.~\ref{fig:evryADD0} and Fig.~\ref{fig:evryCompBDD}, we can see that without adding any additional vehicles, our method significantly recovers accessibility, in particular blue and green tiles in Fig.~\ref{fig:evryADD0} in regions with the largest losses in accessibility due to disruption (red tiles in Fig.~\ref{fig:evryCompBDD}). Similar findings apply to Choisy-le-Roi (comparing Fig.~\ref{fig:clrADD0} and~\ref{fig:clrCompBDD})





Fig.~\ref{fig:accCompCombined} further corroborates the conclusion; in Fig.~\ref{fig:evryACCAffected} and Fig.~\ref{fig:clrACCAffected}, we first compute the loss of accessibility of each centroid due to disruption. We then sort them from the most affected to the least on the $x$-axis. For each of these centroids, we represent in the $y$-axis the improvement  of accessibility of our method versus the conventional replacement bus, i.e.:
\begin{align}
&&\frac{acc(c_i,\pazocal G^\text{OURS})-acc(c_i,\pazocal G^\text{REPL}) }
{acc(c_i,\pazocal G^\text{REPL})},
&& \forall c_i\in\pazocal C    
\label{eq:REPL-OURS}
\end{align}
\begin{figure*}[!ht]
    \centering
    \begin{subfigure}{0.47\linewidth}
        \centering
        \includegraphics[width=\linewidth]{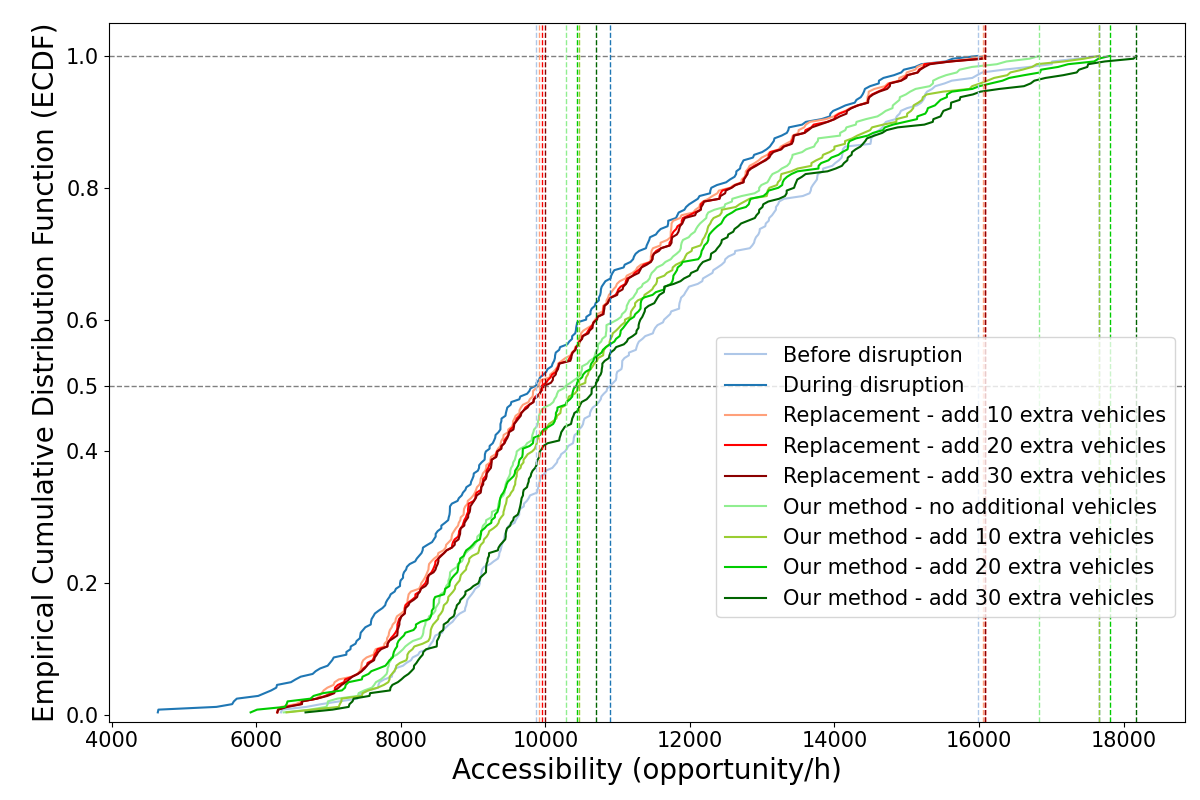}
        \caption{Évry-Courcouronnes.}
        \label{fig::accCompECDF_EVRY}
    \end{subfigure}
    \begin{subfigure}{0.47\linewidth}
        \centering
        \includegraphics[width=\linewidth]{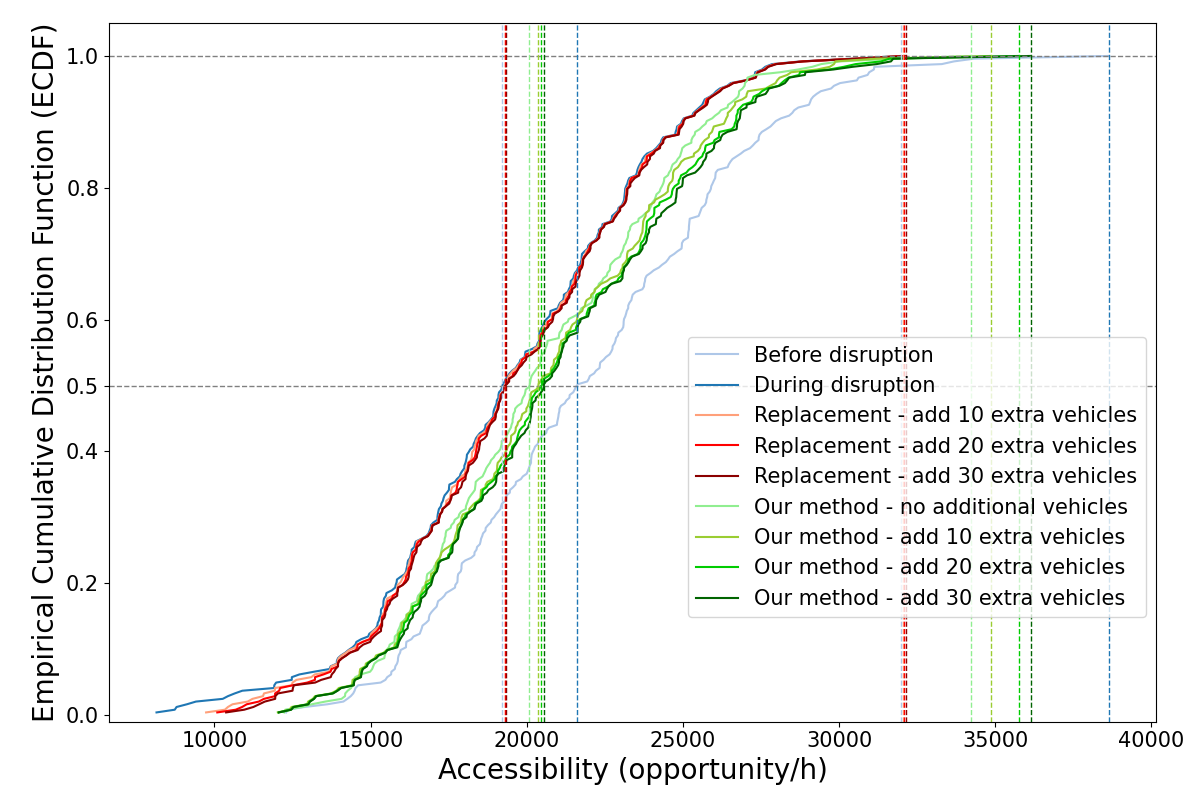}
        \caption{Choisy-le-Roi.}
        \label{fig::accCompECDF_CLR}
    \end{subfigure}
    \caption{Comparison of Cumulative Distribution Functions for Different Scenarios.}
    \label{fig::accCompECDF}
\end{figure*}
\begin{figure*}[!ht]
    \centering
    \begin{subfigure}{0.47\linewidth}
        \centering
        \includegraphics[width=\linewidth]{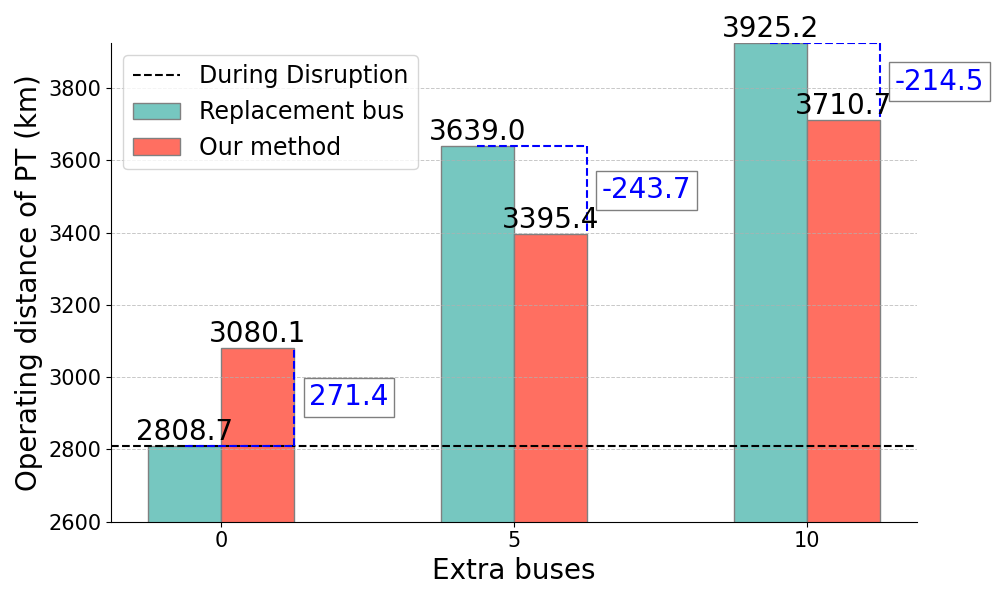}
        \caption{Évry-Courcouronnes}
        \label{fig:evryDistance}
    \end{subfigure}
    \hfill
    \begin{subfigure}{0.47\linewidth}
        \centering
        \includegraphics[width=\linewidth]{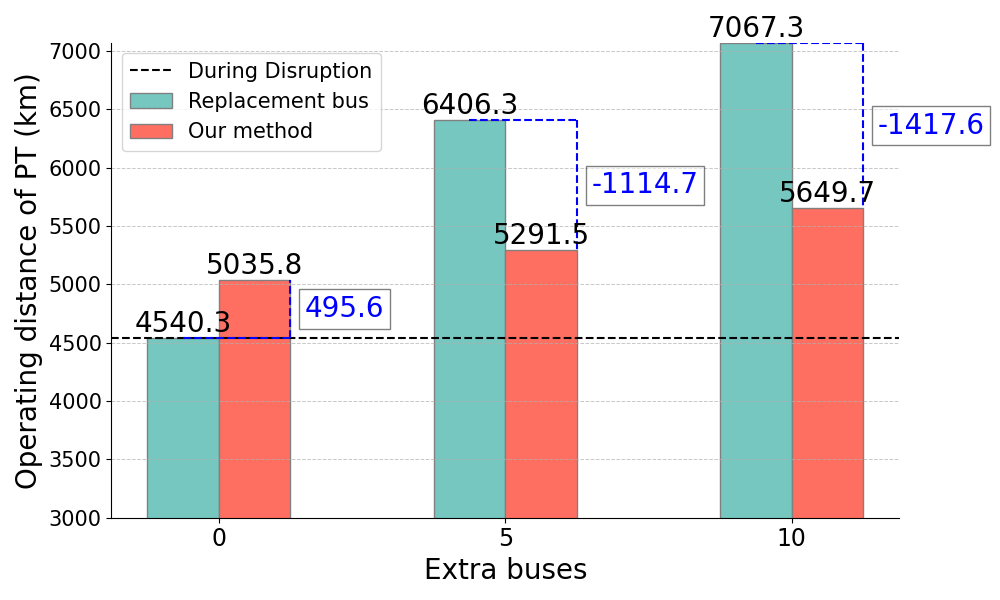}
        \caption{Choisy-le-Roi}
        \label{fig:clrDistance}
    \end{subfigure}
    \caption{Comparison of operating distance (km/h) in two cases}
    \label{fig:distanceComp}
\end{figure*}

we calculate this improvement, by assuming we add a certain number of extra buses (we test~$0, 5$ and~$10$ extra buses).
Note that, when the number of extra buses is $0$, the \emph{Conventional Replacement Bus} scenario corresponds exactly to the \emph{Before Disruption} scenario (indeed, running a bus line that replaces the rail line with $0$ buses operating on it is equivalent to operate no remediation at all). This suggests that our approach improves the accessibility of the vast majority of affected centroids, especially those that are heavily impacted (on the left of the $x$-axis). Even without adding any extra buses, the improvement of our solution over the classic approach goes up to 60\% in the most affected centroids. 
{Figure 5 explicitly shows that compared to the traditional baseline (conventional replacement bus), our method reduces accessibility degradation impact for each centroid (with few exceptions), confirming that while travel time increases may occur, their negative impact on accessibility is mitigated.}
Fig.~\ref{fig:avgACCevry} and Fig.~\ref{fig:avgACCclr} show the accessibility recovered by the Conventional Replacement Bus method and by our solution. In particular, let~$\overline{acc}(\pazocal G)=\frac{1}{|\pazocal C|}\cdot \sum_{c_i\in\pazocal C} acc(c_i, \pazocal G)$ be the average accessibility of PT configuration~$\pazocal G$. We compute in Fig.~\ref{fig::avgACC} the following quantities: $\overline{acc}(\pazocal G^\text{REPL})/\overline{acc}(\pazocal G^\text{O})$ and $\overline{acc}(\pazocal G^\text{OURS})/\overline{acc}(\pazocal G^\text{O})$.
For both Évry-Courcouronnes and Choisy-le-Roi cases, our method is consistently more efficient than the Conventional Replacement Bus method, even when we compare our solution with 0 extra buses with the Conventional Replacement Bus method with 10 extra buses.




Fig.~\ref{fig::accCompECDF} shows the cumulative distribution functions (ECDF) of centroids' accessibility across different scenarios. {The green curves represent our method, and the red curves represent the conventional replacement bus strategy.} We observed that our method substantially increases accessibility for all numbers of extra buses (10, 20, and 30 vehicles). 
For the conventional replacement buses, we can see that even if 30 additional buses are added, it is still not enough to recover the loss in accessibility. Our solution outperforms its results, even with $0$ for additional buses. Moreover, adding extra buses in the conventional replacement bus start only brings negligible improvement, while our solution can consistently increase accessibility when increasing the fleet. Overall, our solution is more effective in narrowing the gap between before-disruption and after-remediation accessibility distribution.

Finally, we study the operational distance covered by the entire fleet of buses in the unit of time. As concerns the Conventional Replacement Bus strategy, we include in the calculation the distance covered by all the regular lines plus the distance covered by the replacement buses. As concerns our solution, instead, we include in the calculation the distance covered by all the regular lines plus the distance covered by the extended lines. We obtained real-world travel distances using Open-Street Map, assuming a set of lines $\pazocal{L}$,  Total travel distance per unit of time = $ \sum_{l\in\pazocal{L}}$ Circle length of line $l \times $ Number of service runs on line $\ell$ per unit of time.
As shown in Fig.~\ref{fig:distanceComp}, our method significantly reduces the total operational distance in both the Évry-Courcouronnes and Choisy-le-Roi case studies. This suggests that our method reduces operational cost~\cite{wang2023integrated} and environmental impact. 

To summarise, our method effectively recovers accessibility loss due to disruption and keeps operating costs lower {, more effectively than the usual replacement bus adopted in practice. Observe that a usual replacement bus is ``easier to intepret'' by the users, as it just follows the same path of the disrupted URT line. A journey performed in the original network could also be replicated when a replacement bus operates instead of the URT line. However, such journey would be so degraded (replacement buses usually have lower frequencies and lower speed, as they are affect by road congestion), that users often prefer to find alternative journeys. Since in our network we re-arrange the PT structure, users of the disrupted line are now forced to calculate a new journey, which may decrease the ``friendliness'' of the system. We believe however, that the loss in friendliness is largely compensated by the much higher efficiency of our approach. Moreover, we take the reasonable assumption that disrupted travelers would not need to carry the burden of recalculating themselves the journey in the re-arranged network: travelers would rather use a smartphone application, which is becoming more and more common nowadays. }

{While our optimized bus configuration demonstrates superior accessibility restoration compared to conventional methods, it serves as a temporary adaptation strategy. The rail service should ultimately be restored due to its inherent advantages in passenger capacity, dedicated rights-of-way, and environmental benefits. However, the accessibility insights gained could inform future network resilience planning, particularly for areas with frequent maintenance or recurring disruptions.}


\section{Conclusion}
This paper presents a method to redesign bus networks during urban rail disruptions by calculating optimal routes for bus line extensions and allocating vehicles among lines. Unlike conventional bus replacement or bridging approaches that focus on operational metrics, our method specifically aims to restore accessibility loss induced by the disruption.

We formulate the bus network redesign problem as an integer program. We then propose a two-stage heuristic resolution method that balances the improvement of accessibility and the kilometers traveled by the fleet. The case studies in Évry-Courcouronnes and Choisy-le-Roi, France, show that our approach is superior to conventional replacement methods, restoring more accessibility with less driving distance, even without additional vehicles.

\section{CONTRIBUTIONS OF THE AUTHORS}
All authors contributed to the concept. Z.G. and A.A. conceptualized the research objectives and developed the methodological framework. Z.G. formulated the optimization model, designed the resolution method, and conducted the analysis. A.A. refined the formulation of the problem and the resolution method.
A.A. and M.E. provided advice and scientific guidance throughout the process of model development, algorithm design, and analysis. A.A. provided funding. All authors reviewed the manuscript.

\section{Acknowledgements}
This work has been supported by the French ANR research project MuTAS (ANR-21-CE22-0025-01).\\
We are grateful to the IBM Academic Initiative for providing IBM ILOG CPLEX Optimization Studio, which was essential for our research. 

\bibliographystyle{unsrt}  
\bibliography{ref}

\end{document}